\documentclass[a4paper,12pt]{article}
\usepackage{amsmath,latexsym}
\usepackage{amsthm}
\usepackage{amsfonts}
\usepackage{dsfont}
\usepackage{upgreek}
\usepackage{tipa}
\usepackage{amssymb}
\usepackage{graphicx}
\usepackage[square,sort,comma,numbers]{natbib}
\usepackage[utf8]{inputenc}
\usepackage{float}
\usepackage{enumitem}   
\newtheorem{theorem}{Theorem}[section]

\newtheorem{definition}[theorem]{Definition}

\theoremstyle{remark}
\newtheorem{remark}[theorem]{Remark}
\vspace{5cm}
\begin{document}

	\begin{center}
		{ 
			{\Large \textbf {Bicomplex Mittag-Leffler Function and Properties    
				}
			}
			\\

			\medskip

			{\sc  Ritu Agarwal$^{1},$\ \ Urvashi Purohit Sharma$^{2}$\ Ravi P. Agarwal$^{3}$}\\
			{\footnotesize $^{1,2}$Department of Mathematics,
				Malaviya National Institute of Technology, Jaipur-302017, INDIA}\\
			$^{3}$Department of Mathematics,
			Texas A{\&}M University - Kingsville 700 University Blvd. Kingsville\\
		}
		{\footnotesize E-mail: {\it $^{1}$ragarwal.maths@mnit.ac.in, $^{2}$urvashius100@gmail.com, $^{3}$ravi.agarwal@tamuk.edu}}
		
	\end{center}
	
	\thispagestyle{empty}

	\hrulefill
	\begin{abstract}
		\indent 
		With the increasing importance of the   Mittag-Leffler function in the physical applications, these days many researchers are studying various generalizations and extensions of the Mittag-Leffler function. In this paper efforts are made  to define bicomplex extension of the  Mittag-Leffler function and also its analyticity  and region of convergence are  discussed. Various properties of the bicomplex  Mittag-Leffler function including integral representation, recurrence relations, duplication	formula and differential relations are established. 
		\end{abstract}
	\hrulefill  
	

	\section{Introduction}
	
		Bicomplex numbers are being studied for quite a long time and a lot of work has been done
	in this area. Cockle \cite{cockle1848,  cockle1849} introduced \textit{Tessarines} between 1848 and 1850 following which Segre \cite{csegre1892} introduced bicomplex numbers.\\
	 Many properties of the bicomplex numbers have been discovered.	During the last few years developments have
	aimed to achieve different algebraic and geometric properties of bicomplex
	numbers and its applications (see,  e.g. \cite{charak2009, rochon2004, rochon2004b, rochon2006, ronn2001, meluna2012}). In the recent developments,  efforts have been done to extend the integral transforms \cite{ ragarwal2017, ragarwal2017c}, holomorphic and meromorphic functions \cite{charak2009, charak2011, charak2011a}  a number of functions like  Polygamma function \cite{goyal2007}, Hurwitz Zeta function \cite{goyal2006a}, Gamma and Beta functions \cite{goyal2006}, Riemann Zeta function \cite{rochon2004},  bicomplex analysis and Hilbert space \cite{romesh2011, romesh2015, romeshk2014, rochon2010b, rochon2010, rochon2010a}  in the bicomplex variable from their complex counterparts. \\
	Recently many researchers have worked on various generalizations and the extensions of the Mittag- Leffler function \cite{gorenflo2014,arshad2018,andric2018,roberto2018,saxena2011}.
\subsection{Bicomplex Numbers}
Segre \cite{csegre1892} defined the set of bicomplex numbers as:
\begin{definition}[Bicomplex Number]
		In terms of real components, the set of bicomplex numbers is defined as
		\begin{equation}
		\mathbb{T} =\{\xi: \xi = x_0 +i_1 x_1 +i_2 x_2 + j x_3~| ~x_0,~x_1,~x_2,~x_3~ \in \mathbb{R} \}, 	
		\end{equation}
		and in terms of complex numbers it can be written as
		\begin{equation}\label{eq:bc}
		\mathbb{T} = \{\xi: \xi =z_1 + i_2z_2~ |~z_1,~z_2 ~\in \mathbb{C}\}.	
		\end{equation} 
	\end{definition}
We shall use the  notations, 
$x_0=\operatorname{Re}(\xi),~x_1=\operatorname{Im}_{i_1}(\xi),~x_2=\operatorname{Im}_{i_2}(\xi),~x_3=\operatorname{Im}_{j}(\xi).$ 

Segre discussed the presence of zero divisors 
 which he called \textit{Nullifics}. He noticed    that the zero divisors in bicomplex numbers constitute two ideals which he called infinite set of nullifics.
	The set of all zero divisors is called  null cone \cite{rochon2004a} defined as follows:
	 
	\begin{equation}\label{eq:s nullcone}
	\mathbb{NC} = \mathbb{O}_2= \{z_1 +z_2i_2 ~|~z_1^2 +z_2 ^2 = 0 \}. 
	\end{equation}

Two non trivial idempotent zero divisors in $\mathbb{T},$ denoted by $e_1$ and $e_2$ and defined as follows \cite{price1991}: \\
$ e_1+ e_2 = 1, ~e_1.e_2   = 0,$	$	e_1 =  \dfrac{1 + i_1i_2}{2}=\dfrac{1+j}{2},~ e_2   =\dfrac  {1 - i_1i_2}{2} =\dfrac{1-j}{2},$   and $e_1^2=e_1,~e_2^2=e_2.$
	\begin{definition}[Idempotent Representation]
		Every element  $\xi\in \mathbb{T}$ has unique idempotent representation defined by 
		\begin{equation}
		 \xi = z_1 + i_2 z_2   = (z_1 - i_1 z_2) e_1 + (z_1 + i_1 z_2) e_2= \xi_1 e_1 + \xi_2 e_2
		\end{equation}  where  $\xi_1 = (z_1 - i_1 z_2)$ and $\xi_2 = (z_1 + i_1 z_2).$ 
	\end{definition} 
Projection mappings  $P_1 : 	\mathbb{T} \rightarrow T_1 \subseteq \mathbb{C},$ $ P_2 : \mathbb{T} \rightarrow T_2 \subseteq \mathbb{C} $ for a bicomplex number $\xi=z_1+i_2z_2$ are defined as  (see, e.g. \cite{rochon2004b}):\\
	\begin{equation}\label{eq:p1}
	P_1(\xi)=	P_1(z_1 + i_2z_2) = P_1[(z_1 - i_1z_2)e_1 + (z_1 + i_1z_2)e_2] = (z_1 - i_1z_2) \in T_1
	\end{equation}
	and
	\begin{equation}\label{eq:p2}
	P_2(\xi)=P_2(z_1 + i_2z_2) = P_2[(z_1 -i_1z_2)e_1 + (z_1 + i_1z_2)e_2] = (z_1 + i_1z_2) \in T_2,
	\end{equation}
	where
	\begin{equation}\label{eq:  a space A1}
	T_1= \{\xi_1= z_1-i_1z_2 \hspace{1mm}| z_1,z_2 \in \mathbb{C}\}~\text{and}~ 	T_2= \{\xi_2= z_1+i_1z_2\hspace{1mm} | z_1,z_2 \in \mathbb{C}\}.
	\end{equation}

\begin{remark}
	The bicomplex space $\mathbb{T}$ can be written as the product 
	\begin{equation}\label{T cart}
	\mathbb{T}=T_1\times_eT_2=\left\lbrace \xi_1e_1+\xi_2e_2,~\xi_1\in T_1,~\xi_2\in T_2\right\rbrace. 
	\end{equation}
\end{remark}
	\begin{definition}[Modulus]
			Let $\xi= z_1+i_2z_2\in\mathbb{T}$ (see, e.g. \cite{rochon2004b}). \\
	The real modulus of $\xi$ is defined as
	\begin{equation}
	|\xi|=\sqrt{|z_1|^2+|z_2|^2}.
	\end{equation}
	The $i_1$-modulus of $\xi$ is defined as
	\begin{equation}
	|\xi|_{i_1}=\sqrt{z_1^2+z_2^2}.
	\end{equation}
	The $i_2$-modulus of $\xi$ is defined as
	\begin{equation}
	|\xi|_{i_2}=\sqrt{(|z_1|^2-|z_2|^2)+2\operatorname{Re}(z_1\bar{z_1})i_2}.
	\end{equation}
	The $j$-modulus of $\xi$ is defined as
	\begin{equation}\label{eq:jmod}
	|\xi|_j=|z_1-i_1z_2|e_1+|z_1+i_1z_2|e_2.
	\end{equation}
\end{definition}

\begin{definition}[Norm]
		Let $\xi= z_1+i_2z_2=\xi_1e_1+\xi_2e_2=x_0+x_1i_1+x_2i_2+x_3j\in\mathbb{T}$ (see, e.g. \cite{meluna2013,riley1953}) then 
	the norm of $\xi$ is given by
	\begin{equation}
	\lVert\xi\rVert=\sqrt{|z_1|^2 + |z_2|^2}=\frac{1}{\sqrt{2}}\sqrt{|\xi_1|^2+|\xi_2|^2}=\sqrt{x_0^2 + x_1^2 + x_2^2 + x_3^2 }.
	\end{equation}
	
\end{definition}
\begin{definition}
		Let $\xi= z_1+i_2z_2=\xi_1e_1+\xi_2e_2=x_0+x_1i_1+x_2i_2+x_3j\in\mathbb{T}$ (see, e.g. \cite{riley1953})  then absolute value of $\xi$ is  denoted by $|\xi|_{abs}$, and is defined as \\
		\begin{equation}
		|\xi|_{abs} = \sqrt{|z_1^2 + z_2^2|} = \sqrt{|(z_1 - i_1 z_2)(z_1 + i_1 z_2)|} = 
		\sqrt{|\xi_1 \xi_2|} = \sqrt{|\xi_1| |\xi_2|}.
		\end{equation}
		
\end{definition}
\begin{definition}[Argument]
	Let $\xi= z_1+i_2z_2=\xi_1e_1+\xi_2e_2=x_0+x_1i_1+x_2i_2+x_3j\in\mathbb{T}$
	then hyperbolic argument (see, e.g. \cite{meluna2014}) of $\xi$ is given by:	\begin{equation}\label{eq:harg}
	\arg_j (\xi)=\arg (\xi_1) e_1+\arg(\xi_1)e_2. 
	\end{equation}
\end{definition}	

Let $U$ be an open set, and $ g : U \subseteq \mathbb{T}\rightarrow \mathbb{T} $ (see, e.g. \cite{rochon2004, ronn2001}). Also
$g (z_1 + i_2z_2) = g_1 (z_1,	z_2) + i_2 g_2 (z_1, z_2)$, then $g$ is $\mathbb{T}$-holomorphic iff $g_1$ and $g_2$ are holomorphic in $U$ and
\begin{equation}\label{eq: CR equation}
\frac{\partial g_1}{\partial z_1} = \frac{\partial g_2}{\partial z_2} \hspace{5mm}\text{and} \hspace{5mm}
\frac{\partial g_2}{\partial z_1} = - \frac{\partial g_1}{\partial z_2} \hspace{2mm}\text{on} \hspace{2 mm}U.
\end{equation}

These equations are called  the bicomplex Cauchy-Riemann equations (abbr. bicomplex CR-equations).
\begin{equation}\label{eq:f'(z)}
g' = \frac{\partial g_1}{\partial z_1} + i_2 \frac{\partial g_2}{\partial z_1}.	 \end{equation}
In the following theorem, Riley \cite{riley1953} studied  the  convergence of bicomplex power series.
\begin{theorem}\label{th:radius of con}
	Let\begin{equation}\label{eq:radius of con}
	N(\xi)=\sqrt{\lVert\xi\rVert^2+\sqrt{\lVert\xi\rVert^4-|\xi|_{abs}^4}}=\max (|\xi_1|,|\xi_2|)
	\end{equation}
	then $N(\xi)$ is a norm and if $\sum_{n=0}^{\infty}a_n\xi^n$ is a power series with component series $\sum_{n=0}^{\infty}b_n\xi_1^n$ and $\sum_{n=0}^{\infty}c_n\xi_2^n,$  $a_n=b_ne_1+c_ne_2$ both have same radius of convergence $R>0$ then $\sum_{n=0}^{\infty}a_n\xi^n$ converges for $N(\xi)<R$ and diverges for $N(\xi)>R,$ where $\lVert\xi\rVert=\frac{1}{\sqrt{2}}\sqrt{|\xi_1|^2+|\xi_2|^2}$ and $|\xi|_{abs}=\sqrt{|\xi_1||\xi_2|}.$
\end{theorem}
In the following theorem, Ringleb \cite{ringlab1933} (see, e.g. \cite{riley1953})  discussed the anlyticity of a bicomplex	function w.r.t. its idempotent complex component functions. This theorem plays a vital role while discussing the convergence of the bicomplex functions.
\begin{theorem}[Decomposition theorem of Ringleb \cite{ringlab1933}]\label{th: ringleb}
	Let $f(z)$ be analytic in a region $\mathbb{T},$ and let $ T_1$ and $T_2$ be the component regions of $\mathbb{T},$ in the $\xi_1$ and $\xi_2$ planes, respectively. Then there exists a unique pair of complex-valued analytic functions, $f_1(\xi_1)$ and $f_2(\xi_2)$, defined in $ T_1$ and $T_2$, respectively, such that \\
	\begin{equation}\label{u1}
	f(z) = f_1(\xi_1)e_1 + f_2(\xi_2) e_2
	\end{equation}
	for all $\xi$ in $\mathbb{T}.$ Conversely, if $f_1(\xi_1)$ is any complex-valued analytic function in
	a region $T_1$ and $ f_2(\xi_2)$ any complex-valued analytic function in a region $T_2$ then the bicomplex-valued function f(z) defined by the equation (\ref{u1})  is an analytic function of the bicomplex variable $\xi$ in the product-region $\mathbb{T}= T_1\times_eT_2$ (equation (\ref{T cart})).
\end{theorem}
In the Theorem \ref{th:int price}, Price \cite{price1991} studied the integration in bicomplex domain w.r.t.  its idempotent representaion. This theorem plays a basic role in the study of integrals of the bicomplex function.
\begin{theorem}\label{th:int price}
	Let $X\subseteq \mathbb{T}.$
	Let $C_1,~C_2$ be two curves defined as 
	\begin{equation}
	C_1:z_1-i_2z_2=z_1(t)-i_1z_2(t)=\xi_1=\xi_1(t),~a\le t\le b.
	\end{equation}	
	\begin{equation}
	C_2:z_1+i_2z_2=z_1(t)+i_1z_2(t)=\xi_2=\xi_2(t),a\le t\le b.
	\end{equation}		
	which have continuous derivatives and whose traces are in $X_1\subseteq T_1,~X_2\subseteq T_2$ respectively and let $C$ be the curve with trace in $X$ which is defined as 
\begin{equation}
C:\xi(t)=\xi_1(t)e_1+\xi_2(t)e_2,~a\le t\le b
\end{equation}	
Then the integrals of $f,~f_1,~f_2$ on the curves $C,~C_1$ and $C_2$ exists and 
\begin{equation}
\int_{C}f(\xi)d\xi=\int_{C_1}f_1(z_1-i_1z_2)d(z_1-i_1z_2)e_1+\int_{C_2}f_1(z_1+i_1z_2)d(z_1+i_1z_2)e_2
\end{equation}
or 
\begin{equation}
\int_{C}f(\xi)d\xi=\int_{C_1}f_1(\xi_1)d(\xi_1)e_1+\int_{C_2}f_2(\xi_2)d(\xi_2)e_2.
\end{equation}

\end{theorem}

We would require the definition of the bicomplex gamma function
	  defined by  Goyal et al. \cite{goyal2006},   in the Euler product form as follows:
	\begin{equation}\label{eq gamma}
	\frac{1}{\Gamma \xi}= \xi e^{\gamma\xi}\prod_{n=1}^{\infty}\left( \left( 1+\frac{\xi}{n}\right) exp \left( -\frac{\xi}{n}\right) \right),~\xi\in\mathbb{T}
	\end{equation}
	 provided that $z_1 \neq \frac{-(m+l)}{2},$ and $z_2 \neq i_1(\frac{l-m}{2})$ where $m,~l \in \mathbb{N} \cup \{0\}.$
	 The Euler	constant	$\gamma(0 \leq \gamma \leq 1)$ is given by
	
	\begin{equation}
	\gamma= \lim\limits_{n \to \infty} (H_n-\log n),~	H_n= \sum_{k=1}^{n}\frac{1}{k}.
	\end{equation}	
	Also, in idempotent form 
	\begin{equation}
	\Gamma \xi=\Gamma \xi_1 e_1+\Gamma \xi_2e_2,~\xi\in\mathbb{T},
	\end{equation}
and in the integral form (see, e.g.\cite{goyal2006}), for $p=p_1e_1+p_2e_2,~ p_1,p_2\in\mathbb{R}^+ .$
\begin{equation}\label{eq:gamma}
\Gamma \xi= \int_{H}e^{-p}p^{\xi-1}dp=\left(\int_{0}^{\infty}e^{-p_1}p_1^{\xi_1-1}dp_1\right) e_1+\left(\int_{0}^{\infty}e^{-p_2}p_2^{\xi_2-1}dp_2\right) e_2
\end{equation}
where $H=(\gamma_1,\gamma_2)$ and $\gamma_1:0 ~\text{to}~ \infty,~\gamma_2:0 ~\text{to}~ \infty.$\\
\noindent\textbf{{Mittag-Leffler Function  and its Properties}}\\
The Mittag-Leffler function (M-L function) comes intrinsically in the study  of the fractional calculus.   The importance of the M-L function in science and engineering is continuously increasing. It is very useful in the area of fractional modeling.


The one parameter  M-L function  defined by Mittag-Leffler \cite{mittag1905} is given by

\begin{equation}\label{eq 2}
\mathbb{E}_{\upalpha}(z)=\sum_{k=0}^{\infty}\frac{z^k}{\Gamma(\upalpha k+1 )}  ,~ \operatorname{Re}(\upalpha)>0 ,~~z \in \mathbb{C}.
\end{equation}
It comes from the Cauchy inequality for the Taylor coefficients and  properties of the Gamma function (see, e.g.\cite[p.18]{gorenflo2014}) that $\exists$ a number $k \ge 0$ and a positive number $r(k)$ such that
\begin{equation}\label{eq: maxE}
M_{\mathbb{E}_\upalpha(r)}=\max_{|z|=r}|\mathbb{E}_\upalpha(z)|<e^{r^k},~\forall r>r(k).
\end{equation}
hence $\mathbb{E}_\upalpha(z)$ is an entire function of finite order.\\
For each $\operatorname{Re}(\upalpha) > 0$ the order $\rho$ and type $\sigma$ of M-L function (\ref{eq 2}) is given by\\

\begin{equation}\label{eq: order}
\rho= \lim \sup_{k \to \infty}\frac{k \log k}{\log\frac{1}{|a_k|}}=\frac{1}{\operatorname{Re}(\upalpha)}
\end{equation}
and 
\begin{equation}\label{eq:type}
\sigma = \frac{1}{e\rho}\lim \sup_{k \to \infty}(k|a_k|^{\frac{\rho}{k}})=1.
\end{equation}


\section{Bicomplex one-parameter  Mittag-Leffler function}
Here, we introduce the bicomplex one parameter Mittag-Leffler function defined by 
	\begin{equation}\label{eq: bc m2}
	\mathbb{E}_{\upalpha}(\xi)=\sum_{k=0}^{\infty}\frac{\xi^k}{\Gamma(\upalpha k+1 )},  \end{equation}
where $\xi, \upalpha\in \mathbb{T},~\xi= z_1+i_2z_2$ and $|\operatorname{Im_j}(\upalpha)|<\operatorname{Re}(\upalpha).$ 


	The definition of bicomplex M-L function is well justified by the following theorem.
	\begin{theorem}
		Let $\xi, \upalpha\in \mathbb{T}$ where $\xi= z_1+i_2z_2= \xi_1e_1+\xi_2e_2$, $\upalpha= \upalpha_1e_1+\upalpha_2e_2= a_0+i_1a_1+ i_2a_2+i_1i_2a_3,$ with $|\operatorname{Im_j}(\upalpha)|<\operatorname{Re}(\upalpha)$.
		Then
		\begin{equation}\label{eq: bc m2}
		\mathbb{E}_{\upalpha}(\xi)=\sum_{k=0}^{\infty}\frac{\xi^k}{\Gamma(\upalpha k+1 )}.  \end{equation}
	
	\end{theorem}
	\begin{proof}
		Consider the function
		\begin{equation}\label{eq: ml}
		\mathbb{E}_{\upalpha}(\xi)=\sum_{k=0}^{\infty}\frac{\xi^k}{\Gamma(\upalpha k+1 )}.  \end{equation}
		By using the idempotent representation
			\begin{equation}\label{eq: bc ml one para}
			\begin{split}
				\mathbb{E}_{\upalpha}(\xi)&=\sum_{k=0}^{\infty}\frac{\xi_1^k}{\Gamma(\upalpha_1 k+1 )}e_1+\sum_{k=0}^{\infty}\frac{\xi_2^k}{\Gamma(\upalpha_2 k+1 )}e_2\\
				&=\mathbb{E}_{\upalpha_1}(\xi_1)e_1+\mathbb{E}_{\upalpha_2}(\xi_2)e_2,
			\end{split}
	\end{equation}
		where $\xi=\xi_1 e_1+\xi_2e_2$ and $\upalpha=\upalpha_1e_1+\upalpha_2 e_2.$
		
Now,		
		 
		\begin{equation}\label{eq cml1}
	\mathbb{E}_{\upalpha_1}(\xi_1)=	\sum_{k=0}^{\infty}\frac{\xi_1^k}{\Gamma(\upalpha_1 k+1 )},
		\end{equation}
		is the  complex M-L function  convergent for $\operatorname{Re}(\upalpha_1)>0 ,~~\xi_1 \in \mathbb{C}.$\\
		Similarly,
		\begin{equation}\label{eq cml1}
     \mathbb{E}_{\upalpha_2}(\xi_2)=		\sum_{k=0}^{\infty}\frac{\xi_2^k}{\Gamma(\upalpha_2 k+1 )},
		\end{equation}
		is also complex M-L function  convergent for $\operatorname{Re}(\upalpha_2)>0 ,~~\xi_2 \in \mathbb{C}.$\\
		Since $ \mathbb{E}_{\upalpha_1}(\xi_1) $ and $\mathbb{E}_{\upalpha_2}(\xi_2)$ are convergent in $T_1,~T_2$ respectively, by Ringleb decomposition theorem  (\ref{eq: ml}) is also convergent in $\mathbb{T}.$\\
		Further, Let
\begin{equation}
\begin{split}
\upalpha&= a_0+i_1a_1+ i_2a_2+i_1i_2a_3\\
&=\upalpha_1e_1+\upalpha_2e_2,
\end{split}
\end{equation}
where $\upalpha_1= (a_0+a_3) +i_1(a_1-a_2)$ and $\upalpha_2= (a_0-a_3) +i_1(a_1+a_2).$

 Since
 	$\operatorname{Re}(\upalpha_1)>0$ and  $\operatorname{Re}(\upalpha_2)  >0$
 	\begin{eqnarray}
 	&\Rrightarrow & a_0+a_3 >0 ~\text{and}~ a_0-a_3 >0.\\
 	&	\Rrightarrow& ~ |a_3|<a_0.\label{eq:re(alpha)g0}\\
 	&\Rrightarrow&|\operatorname{Im_j}(\upalpha)|<\operatorname{Re}(\upalpha).
 	\end{eqnarray}
	
This completes the proof.		
\end{proof}
	By substituting the value of the bicomplex gamma function defined by  equation (\ref{eq gamma}) in the equation (\ref{eq: bc ml one para}) we get the following representation for Mittag-Leffler function:

	\begin{theorem}
	Let $\xi, \upalpha\in \mathbb{T}$ where $\xi= z_1+i_2z_2= \xi_1e_1+\xi_2e_2$, $\upalpha= a_0+i_1a_1+ i_2a_2+i_1i_2a_3=\upalpha_1e_1+\upalpha_2e_2,$ with $|\operatorname{Im_j}(\upalpha)|<\operatorname{Re}(\upalpha).$
	then
\begin{equation}
	\mathbb{E}_{\upalpha}(\xi)=\sum_{k=0}^{\infty}{\xi^k} (\upalpha k+1) e^{\gamma(\upalpha k+1 )}\prod_{n=1}^{\infty}\left( \left( 1+\frac{(\upalpha k+1 )}{n}\right) exp \left( -\frac{(\upalpha k+1 )}{n}\right) \right).  
	\end{equation}
	\
\end{theorem}
\begin{remark}
Also, in integral form, bicomplex M-L function can be represented as 
\begin{equation}
\mathbb{E}_{\upalpha}(\xi)= \sum_{k=0}^{\infty}\frac{\xi^k}{\displaystyle\int_{H}e^{-p}p^{\upalpha k }dp}.  
\end{equation}

where $H=(\gamma_1,\gamma_2)$ as defined in (\ref{eq gamma}).\\
	
\end{remark}

For different values of the  $\upalpha$ we obtain  various bicomplex functions as special cases. To mention, a few are:
\begin{enumerate}
	\item For $\upalpha=0$ we get bicomplex binomial function
	$	\mathbb{E}_{0}(\xi)=\frac{1}{1-\xi},~\lVert\xi\rVert<1.$
	\item For $\upalpha=1$ we get bicomplex exponential function
	$	\mathbb{E}_{1}(\pm\xi)=e^{\pm\xi}.$
	\item For $\upalpha=2$ we get  bicomplex cosine function
	$\mathbb{E}_{2}(-\xi^2)=\cos\xi.$
	\item For $\upalpha=2$ we get  bicomplex hyperbolic cosine  function
	$\mathbb{E}_{2}(\xi^2)=\cosh\xi.$	
	\item For $\upalpha=3$ we get following function\\
	$\mathbb{E}_{3}(\xi)=\frac{1}{2}\left(e^{\xi^{1/3}}+2e^{-(1/2)\xi^{1/3}}\cos\left(\frac{\sqrt{3}}{2}\xi^{1/3} \right)  \right).$	
	\item For $\upalpha=4$ we get following bicomplex relation\\
	$\mathbb{E}_{4}(\xi)=\frac{1}{2}\left( \cos(\xi^{1/4})+\cosh(\xi^{1/4})\right).$	
	\end{enumerate}


\begin{theorem}
	The bicomplex Mittag-Leffler function defined in equation (\ref{eq: bc m2}) satisfies bicomplex Cauchy - Riemann equations.
\end{theorem}

\begin{proof}
	By the result (\ref{eq: bc ml one para}) we have,
	\begin{equation}
	\begin{split}
	\mathbb{E}_\upalpha(\xi)&=\mathbb{E}_{\upalpha_1}(\xi_1)e_1+\mathbb{E}_{\upalpha_2}(\xi_2)e_2\\
	&=\mathbb{E}_{\upalpha_1}(z_1 - i_1z_2)e_1+\mathbb{E}_{\upalpha_2}(z_1 + i_1z_2)e_2\\
	&=\mathbb{E}_{\upalpha_1}(z_1 - i_1z_2)\left( \frac{1 + i_1i_2}{2}\right) +\mathbb{E}_{\upalpha_2}(z_1 + i_1z_2)\left( \frac{1 - i_1i_2}{2}\right) \\
	&=\left( \frac{1}{2}\left( \mathbb{E}_{\upalpha_1}(z_1 - i_1z_2)+\mathbb{E}_{\upalpha_2}(z_1 +  i_1z_2)\right) \right)\\
	&~~+i_2 \left( \frac{i_1}{2}\left( \mathbb{E}_{\upalpha_1}(z_1 - i_1z_2)-\mathbb{E}_{\upalpha_2}(z_1 +  i_1z_2)\right) \right)\\
	&=f_1(z_1,z_2)+i_2f_2(z_1,z_2).\\
	\end{split}
	\end{equation}
	where 
	$f_1(z_1,z_2)= \dfrac{1}{2}\left( \mathbb{E}_{\upalpha_1}(z_1 - i_1z_2)+\mathbb{E}_{\upalpha_2}(z_1 +  i_1z_2)\right),$\\
and 	$f_2(z_1,z_2)=  \dfrac{i_1}{2}\left( \mathbb{E}_{\upalpha_1}(z_1 - i_1z_2)-\mathbb{E}_{\upalpha_2}(z_1 +  i_1z_2)\right).$\\
	 $\mathbb{E}_{\upalpha_i}(i=1,2)$ are complex M-L functions.\\
	Now,
	\begin{eqnarray}
	\frac{\partial f_1}{\partial z_1 }&=&\frac{1}{2}\left( \mathbb{E}_{\upalpha_1}'(z_1 - i_1z_2)+\mathbb{E}_{\upalpha_2}'(z_1 +  i_1z_2)\right), \nonumber\\
	\frac{\partial f_1}{\partial z_2 }&=&\frac{-i_1}{2}\left( \mathbb{E}_{\upalpha_1}'(z_1 - i_1z_2)-\mathbb{E}_{\upalpha_2}'(z_1 +  i_1z_2)\right),\nonumber\\
	\frac{\partial f_2}{\partial z_1 }&=&\frac{i_1}{2}\left( \mathbb{E}_{\upalpha_1}'(z_1 - i_1z_2)-\mathbb{E}_{\upalpha_2}'(z_1 +  i_1z_2)\right),\nonumber\\
	\frac{\partial f_2}{\partial z_2 }&=&\frac{1}{2}\left( \mathbb{E}_{\upalpha_1}'(z_1 - i_1z_2)+\mathbb{E}_{\upalpha_2}'(z_1 +  i_1z_2)\right).\nonumber 
	\end{eqnarray}
	from the above equations it can be observed that
	\begin{equation}
	\frac{\partial f_1}{\partial z_1} = \frac{\partial f_2}{\partial z_2} \hspace{5mm}\text{and} \hspace{5mm}
	\frac{\partial f_2}{\partial z_1} = - \frac{\partial f_1}{\partial z_2}.
	\end{equation}
	Hence, bicomplex Cauchy-Riemann equations are satisfied by the bicomplex M-L function.
\end{proof}

\begin{theorem}
The bicomplex M-L function $\mathbb{E}_{\upalpha}(\xi),$ $|\operatorname{Im_j}(\upalpha)|<\operatorname{Re}(\upalpha)$
 is an entire function in the bicomplex domain.	
\end{theorem}
\begin{proof}	
	Let  $ \sum_{n=0}^{\infty} a_n\xi^n$ represents a bicomplex power series where $a_n,~\xi,\in \mathbb{T},~a_n= b_ne_1+ c_ne_2,$  $\xi=\xi_1e_1+\xi_2e_2  $. Then by Ringleb decomposition theorem \ref{th: ringleb}, the series
	\begin{equation}
	 \sum_{n=0}^{\infty} a_n\xi^n=\left( \sum_{n=0}^{\infty} b_n\xi_1^n\right)  e_1+ \left( \sum_{n=0}^{\infty} c_n\xi_2^n\right) e_2
	\end{equation}
	 converges iff 
	$ \sum_{n=0}^{\infty} b_n\xi_1^n $ and $ \sum_{n=0}^{\infty} c_n\xi_2^n $ converge in the complex domains (see, e.g. \cite{riley1953}).\\	
	Now from equation (\ref{eq: bc ml one para}), the Mittag -Leffler function can be decomposed as 
\begin{equation}
	\mathbb{E}_\upalpha(\xi)=\mathbb{E}_{\upalpha_1}(\xi_1)e_1+\mathbb{E}_{\upalpha_2}(\xi_2)e_2
		\end{equation}
	Since $\mathbb{E}_{\upalpha_1}(\xi_1)=\sum_{k=0}^{\infty}\frac{\xi_1^k}{\Gamma(\upalpha_1 k+1 )},\operatorname{Re}(\upalpha_1)>0$ and $\mathbb{E}_{\upalpha_2}(\xi_2)=\sum_{k=0}^{\infty}\frac{\xi_2^k}{\Gamma(\upalpha_2 k+1 )},~\operatorname{Re}(\upalpha_2)>0$ are  complex Mittag Leffler functions with infinite  radius of covergence  ( say $R$)  \cite[p.18]{gorenflo2014}. Then
			\begin{equation}
			|\xi_1|<R,~|\xi_2|<R.
			\end{equation}
		From equation (\ref{eq:radius of con}),\\
		\begin{equation}
	N(\xi)=\sqrt{\lVert\xi\rVert^2+\sqrt{\lVert\xi\rVert^4-|\xi|_{abs}^4}}=\max (|\xi_1|,|\xi_2|)<R.
			\end{equation}
				Hence from Theorem \ref{th:radius of con}, $\mathbb{E}_{\upalpha}(\xi)$ converges in the bicomplex domain and  has infinite radius of convergence \cite{riley1953}. Since complex M-L function is  entire function in $\mathbb{C}$  the bicomplex M-L function is  an entire function in $\mathbb{T}$   (Riley \cite[p.141]{riley1953}).  \end{proof} 
\begin{theorem}[Order and Type]
	The bicomplex Mittag-Leffler function $\mathbb{E}_\upalpha(\xi), ~ \xi, ~\upalpha \in\mathbb{T}$ is an entire function of finite order $\rho= \frac{a_0-a_3j}{(a_0^2-a_3^2)} $ and type $\sigma=1.$
\end{theorem}
\begin{proof}
	
	 From equation (\ref{eq: bc ml one para})
	\begin{equation}
	\mathbb{E}_\upalpha(\xi)=\mathbb{E}_{\upalpha_1}(\xi_1)e_1+\mathbb{E}_{\upalpha_2}(\xi_2)e_2
	\end{equation}
	Here $\mathbb{E}_{\upalpha_1}(\xi_1)$ and $\mathbb{E}_{\upalpha_2}(\xi_2)$ are the complex Mittag-Leffler functions for $\operatorname{Re}(\upalpha_1)>0,~ \operatorname{Re}(\upalpha_2)>0$ respectively. Since $\mathbb{E}_{\upalpha_1}(\xi_1), ~\mathbb{E}_{\upalpha_2}(\xi_2)$  are entire functions,   there exists  numbers $k_1,~ k_2 \ge 0$ and  positive numbers $r_1(k_1), ~r_2(k_2),$ such that, from equation (\ref{eq: maxE}), we get
	\begin{equation}
	M_{\mathbb{E}_{\upalpha_1}(r_1)}=\max_{|\xi_1|=r_1}|\mathbb{E}_{\upalpha_1}(\xi_1)|<e^{r_1^{k_1}},~\forall r_1>r_1(k_2)
	\end{equation}
	and
	\begin{equation}
	M_{\mathbb{E}_{\upalpha_2}(r_2)}=\max_{|\xi_2|=r_1}|\mathbb{E}_{\upalpha_2}(\xi_2)|<e^{r_2^{k_2}},~\forall r_2>r_2(k_2)
	\end{equation}
	Let $r=\max(r_1,r_2)$ and $k=\max(k_1, k_2)$ then
	\begin{equation}
	M_{\mathbb{E}_{\upalpha_1}(r_1)}=\max_{|\xi_1|=r_1}|\mathbb{E}_{\upalpha_1}(\xi_1)|<e^{r_1^{k_1}}\le e^{r^k},
	\end{equation}
	and
	\begin{equation}
	M_{\mathbb{E}_{\upalpha_2}(r_2)}=\max_{|\xi_2|=r_1}|\mathbb{E}_{\upalpha_2}(\xi_2)|<e^{r_2^{k_2}}\le e^{r^k}
	\end{equation}	
	\begin{equation}
	\begin{split}
	M_{\mathbb{E}_\upalpha(r)}&=\max_{|\xi|_j=r}|\mathbb{E}_{\upalpha}(\xi)|_j~~~~~~~~~~~~~~~~~\text{[J-modulus of bicomplex number]}\\ &=\max_{|\xi_1|=r}|\mathbb{E}_{\upalpha_1}(\xi_1)|e_1+\max_{|\xi_2|=r}|\mathbb{E}_{\upalpha_2}(\xi_2)|e_2\\
	&\le e^{r^{k}}e_1+e^{r^{k}}e_2,~~~\forall r>r(k_)\\
	&=e^{r^k}.
	\
	\end{split}
	\end{equation}
	Hence $\mathbb{E}_\upalpha(z)$ is an entire function of finite order.\\
	
For the bicomplex M-L function  $\xi,~ \upalpha\in \mathbb{T}, ~|\operatorname{Im_j}(\upalpha)|<\operatorname{Re}(\upalpha).$ (from equation (\ref{eq:re(alpha)g0})) the order $\rho$ is given by 
\begin{equation}
\begin{split}
\rho
&=\lim \sup_{k \to \infty}\frac{k \log k}{\log \Gamma(\upalpha k+1)}\\
&=\left( \lim \sup_{k \to \infty}\frac{k \log k}{\log \Gamma(\upalpha_1 k+1)}\right) e_1+\left( \lim \sup_{k \to \infty}\frac{k \log k}{\log \Gamma(\upalpha_2 k+1)}\right) e_2.\\
\end{split}
\end{equation}
Now, from equation (\ref{eq: order}),
\begin{equation}
\begin{split}
\rho&=\left( \frac{1}{\operatorname{Re}(\upalpha_1)}\right) e_1+\left( \frac{1}{\operatorname{Re}(\upalpha_2)}\right) e_2=\left( \frac{1}{a_0+a_3}\right) e_1+\left( \frac{1}{a_0-a_3}\right) e_2
=\frac{a_0-a_3j}{(a_0^2-a_3^2)}.
\end{split}
\end{equation}
~~~~~~~~~~~~~~~~~~~~~~~~~~~~~~~~~~~~~~~~~~~~~~~~~~~~~~~~~~~~~~~~~~~~~~[$\because$ $a_0>|a_3|\Rightarrow a_0^2-a_3^2\ne0$ ]
The type $\sigma$ of the bicomplex M-L function $\mathbb{E}_\upalpha(\xi)$ is given by
\begin{equation}
\begin{split}
\sigma& = \frac{1}{e\rho}\lim \sup_{k \to \infty}(k|a_k|_j^{\frac{\rho}{k}})\\
& = \frac{1}{e\rho}\lim \sup_{k \to \infty}(k\Big{|} \frac{1}{\Gamma (\upalpha k+1)}\Big{|}_j^{\frac{\rho}{k}})\\
&=\left( \frac{1}{e\rho}\lim \sup_{k \to \infty}\left( k\Big{|}\frac{1}{\Gamma (\upalpha_1 k+1)}\Big{|}^{\frac{\rho}{k}}\right) \right) e_1+\left( \frac{1}{e\rho}\lim \sup_{k \to \infty}\left( k\Big{|}\frac{1}{\Gamma (\upalpha_2 k+1)}\Big{|}^{\frac{\rho}{k}}\right) \right) e_2\\
&=1.e_1+1.e_2 ~~~~~~~~~~[\text{using equation (\ref{eq:type})}]\\
&=1.
\end{split}
\end{equation}


\end{proof}
\begin{remark}
	There are different moduli such as real, $i_1,$ $i_2$ and $j$ modulus are defined for a bicomlex number (see, e.g. \cite{rochon2004b}). In this paper,  $j$ modulus has been used for the calculation, since it provides expression in terms of idempotent components of the complex modulus. 
\end{remark}
\subsection{Properties of Bicomplex Mittag-Leffler Function}
Integral representation for the complex M-L  function $\mathbb{E}_{\upalpha}(z)$ is given by (see, e.g.\cite[p.209]{erdelyi1955}):
\begin{equation}\label{eq: int}
\int_{0}^{\infty}e^{-t}\mathbb{E}_\upalpha(t^\upalpha z)dt=\frac{1}{1-z},
~z \in\mathbb{C},~\upalpha\ge0.
\end{equation}
The above integral converges in  the unit circle and is bounded by the line $\operatorname{Re}z^{1/\upalpha}=1.$

\begin{theorem}[Integral Representation for bicomplex M-L function]
	Let $\xi\in \mathbb{T}$ where $\xi= z_1+i_2z_2=\xi_1e_1+\xi_2e_2$ and $\upalpha\ge0,~\lVert\xi\rVert<1$ then
	
	\begin{equation}
	\int_{0}^{\infty}e^{-t}\mathbb{E}_\upalpha(t^\upalpha \xi)dt=\frac{1}{1-\xi}.
	\end{equation}
	The above integral converges in  the unit circle and is bounded by the plane 
	$\operatorname{Re}(\xi^{1/\upalpha})=1,~\operatorname{Im_j}(\xi)=0 .$
	
\end{theorem}
\begin{proof}
	By the integral representation (\ref{eq: int}) and the result (\ref{eq: bc ml one para}) we have	for  $\xi\in \mathbb{T}$ where $\xi= z_1+i_2z_2=\xi_1e_1+\xi_2e_2,~ ~\upalpha\ge0$ and $|\xi_1|<1,~ |\xi_2|<1$ 
	\begin{equation}
	\begin{split}
	\int_{0}^{\infty}e^{-t}\mathbb{E}_\upalpha(t^\upalpha \xi)dt&=\left( 	\int_{0}^{\infty}e^{-t}\mathbb{E}_{\upalpha}(t^{\upalpha} \xi_1)dt\right) e_1+\left( 	\int_{0}^{\infty}e^{-t}\mathbb{E}_{\upalpha}(t^{\upalpha} \xi_2)dt\right) e_2\\
	&=\left( \frac{1}{1-\xi_1}\right) e_1+\left( \frac{1}{1-\xi_2}\right) e_2\\
	&=\frac{1}{1-(\xi_1e_1+\xi_2e_2)}\\
	&=\frac{1}{1-\xi}.
	\end{split}
	\end{equation}
	In terms of real components 
$	\xi=x_0 +i_1 x_1 +i_2 x_2 + j x_3=\xi_1e_1+\xi_2e_2.$\\
		Here $\xi_1=(x_0+x_3)+i_1(x_1-x_2),~\xi_2=(x_0-x_3)+i_1(x_1+x_2).$\\
	Since,
	
	\begin{eqnarray}
	&&	|\xi_1|<1~ \text{and}~ |\xi_2|<1\nonumber\\
	&\Rightarrow&	\sqrt{(x_0+x_3)^2+(x_1-x_2)^2}<1~ \text{and}~\sqrt{(x_0-x_3)^2+(x_1+x_2)^2}<1\nonumber\\
	&\Rightarrow&\sqrt{x_0^2+x_1^2+x_2^2+x_3^2+2x_0x_3-2x_1x_2}<1\nonumber \\  
	& ~~~~~~~~~ &\text{and}~\sqrt{x_0^2+x_1^2+x_2^2+x_3^2-2x_0x_3+2x_1x_2}<1\nonumber \\
	&\Rightarrow&\sqrt{x_0^2+x_1^2+x_2^2+x_2^2}<1\nonumber\\
	&\Rightarrow&\lVert\xi\rVert
	<1.\nonumber
	\end{eqnarray}
	Also,
	\begin{eqnarray}
	&&\operatorname{Re}\xi_1^{1/\upalpha}=1,~\operatorname{Re}\xi_2^{1/\upalpha}=1\nonumber\\
	&\Rightarrow&	(x_0+x_3)^{\frac{1}{\upalpha}}=1,~(x_0-x_3)^{\frac{1}{\upalpha}}=1\nonumber\\
	&\Rightarrow&	(x_0+x_3)=1,~(x_0-x_3)=1\nonumber\\
	&\Rightarrow&x_0=1,~x_3=0\nonumber \\ 
	&\Rightarrow& \operatorname{Re}\xi=1,~\operatorname{Im_j}\xi=0 .\nonumber \\
	& \Rightarrow & \operatorname{Re}\xi^{1/\upalpha}=1,~\operatorname{Im_j}\xi=0 .\nonumber \\
	\end{eqnarray}
	
\end{proof}

The complex Mittag-Leffler function has following integral representation (see, e.g. \cite{erdelyi1955})
\begin{equation}\label{eq:i}
\mathbb{E}_\upalpha(z) =\frac{1}{2\pi i}\int_{\Omega}\frac{t^{\upalpha-1} e^t}{t^{\upalpha}-z}dt,~\upalpha>0,~ z \in\mathbb{C}
\end{equation}
where the path of integration $\Omega$ is a loop starting and ending at $-\infty$ and encircling the circular disk $|t| \le |z|^{1/\upalpha}$	  in the positive sense, $ |\arg t| <\pi$ on $\Omega$. 

\begin{theorem}
	Let $\xi,~\omega\in \mathbb{T}$ where $\xi= z_1+i_2z_2=\xi_1e_1+\xi_2e_2,~\omega=\omega_1e_1+\omega_2e_2$ 	  then bicomplex Mittag-Leffler function has following integral representation
	  \begin{equation}
	  \mathbb{E}_\upalpha(\xi) =\frac{1}{2\pi i_1}\int_{H}\frac{\omega^{\upalpha-1} e^(\omega)}{\omega^{\upalpha}-\xi}d\omega,~\upalpha>0
	  \end{equation}
	  where the path of integration $H= (\Omega_1,\Omega_2)$ and $\Omega_1,\Omega_2$ are  loops starting and ending at $-\infty$ and encircling the circular disks $|\omega_1| \le |\xi_1|^{1/\upalpha},~|\omega_2| \le |\xi_2|^{1/\upalpha},$ respectively,	  in the positive sense, $ |\arg \omega_1| <\pi$ on $\Omega_1$ and $ |\arg \omega_1| <\pi$ on $\Omega_1,~\Omega_2$.
\end{theorem}
\begin{proof}
			By the integral representation (\ref{eq:i}) , result (\ref{eq: bc ml one para}) and  the Theorem \ref{th:int price},  we have	for $\xi,~\omega\in \mathbb{T}$ 
			\begin{equation}
			\begin{split}
			E_{\upalpha}(\xi)&= E_{\upalpha}(\xi_1)e_1+E_{\upalpha}(\xi_2)e_2\\
			&=\frac{1}{2\pi i_1}\int_{\Omega_1}\frac{\omega_1^{\upalpha-1} e^{\omega_1}}{\omega_1^{\upalpha}-\xi_1}d\omega_1e_1+\frac{1}{2\pi i_1}\int_{\Omega_2}\frac{\omega_2^{\upalpha-1} e^{\omega_2}}{{\omega_2}^{\upalpha}-\xi_2}d\omega_2 e_2\\
			&=\frac{1}{2\pi i_1} \int_{(\Omega_1,\Omega_2)} \frac{(\omega_1e_1+\omega_2e_2)^{\upalpha-1} e^{(\omega_1e_1+\omega_2e_2)}}{(\omega_1e_1+\omega_2e_2)^{\upalpha}-(\xi_1e_1+\xi_2e_2)}d(\omega_1e_1+\omega_2e_2)\\
			&=\frac{1}{2\pi i_1}\int_{H}\frac{\omega^{\upalpha-1} e^\omega}{\omega^{\upalpha}-\xi}d\omega.
			\end{split}
			\end{equation}
		The path of integration is $H= (\Omega_1,\Omega_2),$ where $\Omega_1,\Omega_2$ are  loops starting and ending at $-\infty$ and encircling the circular disks $|\omega_1| \le |\xi_1|^{1/\upalpha},~|\omega_2| \le |\xi_2|^{1/\upalpha},$ respectively,	  in the positive sense.\\
	Further, since 
	
		\begin{equation}
		|\arg \omega_1| <\pi ~\text{and}~|\arg \omega_2| <\pi,
	\end{equation}
 from the equations (\ref{eq:jmod}) and (\ref{eq:harg})   we have
\begin{equation}
\arg_j \omega=(\arg \omega_1)e_1+(\arg \omega_2)e_2,
\end{equation}
\begin{equation}
\Rightarrow|\arg_j \omega|_j=|\arg \omega_1|e_1+|\arg \omega_2|e_2<\pi e_1+\pi e_2=\pi.
\end{equation}
	\end{proof}
\noindent Recurrence relation for the complex M-L function $\mathbb{E}_{\upalpha}(z)$ is given by the following  relation where  $p,~q$ are the  relatively prime natural numbers (see, e.g. \cite[p.21]{gorenflo2014})

\begin{equation}\label{eq:rec1}
\mathbb{E}_{p/q}(z)=\frac{1}{q} \sum_{l=0}^{q-1}\mathbb{E}_{1/p}(z^{1/q}e^{\frac{2 \pi li_1}{q}}).
\end{equation}

\begin{theorem}[Recurrence Relation for bicomplex M-L function]
	Let $\xi\in \mathbb{T}$ where $\xi= z_1+i_2z_2 $ and $p,q\in\mathbb{N}$ are relatively prime. Then the bicomplex Mittag-Leffler function satisfies
	\begin{equation}
		\mathbb{E}_{p/q}(\xi)=\frac{1}{q} \sum_{l=0}^{q-1}\mathbb{E}_{1/p}\left( \xi^{1/q}e^{\frac{2 \pi li_1}{q}}\right).
	\end{equation}

\end{theorem}

\begin{proof}
	By the recurrence relation (\ref{eq:rec1}) and the result (\ref{eq: bc ml one para}) we have	for  $\xi= z_1+i_2z_2=\xi_1e_1+\xi_2e_2= (z_1-i_1z_2)e_1+(z_1+i_1z_2)e_2.$
	\begin{equation}
	\begin{split}
	\mathbb{E}_{p/q}(\xi)&=\mathbb{E}_{p/q}(\xi_1)e_1+\mathbb{E}_{p/q}(\xi_2)e_2 ,~ q \in \mathbb{N}\\
	&=\mathbb{E}_{p/q}(z_1-i_1z_2)e_1+\mathbb{E}_{p/q}(z_1+i_1z_2)e_2\\
	&=\left( \frac{1}{q} \sum_{l=0}^{q-1}\mathbb{E}_{1/p}\left( (z_1-i_1z_2)^{1/q}e^{\frac{2 \pi li_1}{q}}\right) \right) e_1+\left(\frac{1}{q} \sum_{l=0}^{q-1}\mathbb{E}_{1/p}\left( (z_1+i_1z_2)^{1/q}e^{\frac{2 \pi li_1}{q}}\right)  \right) e_2 \\
	&=\frac{1}{q} \sum_{l=0}^{q-1}\mathbb{E}_{1/p}\left( \xi^{1/q}e^{\frac{2 \pi li_1}{q}}\right) .
	\end{split}
	\end{equation}
\end{proof}
\noindent Duplication formula for the complex M-L function $\mathbb{E}_{\upalpha}(z)$ is defined as:  (see, e.g. \cite[p.53]{gorenflo2014})

\begin{equation}\label{eq:duplication}
\mathbb{E}_{2\upalpha}(z^2)=\frac{1}{2}\left( \mathbb{E}_\upalpha(z)+\mathbb{E}_\upalpha(-z)\right), ~\operatorname{Re}(\upalpha)>0.
\end{equation}

\begin{theorem}[Duplication Formula for bicomplex M-L function]
	Let $\xi,\upalpha\in \mathbb{T}$ where $\xi= z_1+i_2z_2,~ ~|\operatorname{Im_j}(\upalpha)|<\operatorname{Re}(\upalpha)$ then
\begin{equation}
\mathbb{E}_{2\upalpha}(\xi^2)=\frac{1}{2}\left( \mathbb{E}_\upalpha(\xi)+\mathbb{E}_\upalpha(-\xi)\right).
\end{equation}
\end{theorem}
\begin{proof}
 we have
	for $\xi, \upalpha\in \mathbb{T}$ where $\xi= z_1+i_2z_2=\xi_1e_1+\xi_2e_2= (z_1-i_1z_2)e_ 1+ (z_1+i_1z_2)e_2~and~\upalpha=\upalpha_1e_1+\upalpha_2e_2, ~|\operatorname{Im_j}(\upalpha)|<\operatorname{Re}(\upalpha).$ 
	
	\begin{equation}
	\begin{split}
	\frac{1}{2}\left( \mathbb{E}_\upalpha(\xi)+\mathbb{E}_\upalpha(-\xi)\right)&=\frac{1}{2}\left(	\sum_{k=0}^{\infty}\frac{\xi^k}{\Gamma(\upalpha k+1 )}+	\sum_{k=0}^{\infty}\frac{(-\xi)^k}{\Gamma(\upalpha k+1 )}\right)  \\
	&=\frac{1}{2}	\sum_{k=0}^{\infty}\frac{\left( \xi^k+(-\xi)^k\right) }{\Gamma(\upalpha k+1 )}\\
		&=\frac{1}{2}	\sum_{k=0}^{\infty} \frac{\left( \xi^k(1+(-1)^k)\right) }{\Gamma(\upalpha k+1 )}\\
		&=\frac{1}{2}\left( 2+2\frac{\xi^2}{\Gamma (2\upalpha+1)}+2\frac{\xi^4}{\Gamma (4\upalpha+1)}+2\frac{\xi^6}{\Gamma (6\upalpha+1)}+.........\right) \\
		&=\left( 1+\frac{(\xi^2)^1}{\Gamma (1(2\upalpha)+1)}+\frac{(\xi^2)^2}{\Gamma (2(2\upalpha)+1)}+\frac{(\xi^2)^3}{\Gamma (3(2\upalpha)+1)}+.........\right) \\
		&=\sum_{k=0}^{\infty} \frac{(\xi^2)^k}{\Gamma(2\upalpha k+1 )}\\
	&=	\mathbb{E}_{2\upalpha}(\xi^2).
		\end{split}
	\end{equation}
\end{proof}

\noindent Differential relations   for the complex M-L function $\mathbb{E}_{\upalpha}(z)$ are defined by the following  relations where $p,q\in\mathbb{N}$ are relatively prime (see, e.g. \cite[p.22]{gorenflo2014}):

\begin{equation}\label{eq:diff1}
\left(\frac{d}{dz} \right) ^p \mathbb{E}_p(z^p)=\mathbb{E}_p(z^p),
\end{equation}
\begin{equation}\label{eq:diff2}
\frac{d^p}{dz^p}\mathbb{E}_{p/q}(z^{p/q})=\mathbb{E}_{p/q}(z^{p/q})+\sum_{k=1}^{q-1}\frac{z^{-kp/q}}{\Gamma~(1-kp/q)}.
\end{equation}

\begin{theorem}[Differential Relations for the bicomplex M-L function]
	Let $\xi\in \mathbb{T}$ where $\xi= z_1+i_2z_2$ then for $p,~q,$ relatively prime natural numbers
	\begin{enumerate}[label=(\roman*)]
		\item 
		$\left(\frac{d}{d\xi} \right) ^p \mathbb{E}_p(\xi^p)=\mathbb{E}_p(\xi^p).$
		\item 
		$\frac{d^p}{d\xi^p}\mathbb{E}_{p/q}(\xi^{p/q})= \mathbb{E}_{p/q}(\xi^{p/q})+\sum_{k=1}^{q-1}\frac{\xi^{-kp/q}}{\Gamma~(1-kp/q)}.
		$
	\end{enumerate}
	
\end{theorem}
\begin{proof}(i)
	\begin{equation}
	\begin{split}
	\left(\frac{d}{d\xi} \right) ^p \mathbb{E}_p(\xi^p)&=\left(\frac{d}{d\xi} \right) ^p\sum_{k=0}^{\infty}\frac{\xi^{pk}}{\Gamma(p k+1 )}~~~~~~~~~~~~~~~~~~~~[\text{From definition (\ref{eq: bc m2})}]\\
	&=\sum_{k=1}^{\infty}\frac{\xi^{pk-p}}{\Gamma(p k-p+1 )}\\
		&=\sum_{k=0}^{\infty}\frac{\xi^{pk}}{\Gamma(p k+1 )}~~~~~~~~~~~~~~~~~~~~~~~~~~~~~~[ \text{Replacing} ~ k\rightarrow k+1]\\
			&=\mathbb{E}_p(\xi^p).
	\end{split}
	\end{equation}
\end{proof}
\begin{proof}(ii)
	Again,
	\begin{equation}
	\begin{split}
	\frac{d^p}{d\xi^p}\mathbb{E}_{p/q}(\xi^{p/q})&=\frac{d^p}{d\xi^p}\sum_{k=0}^{\infty}\frac{\xi^{\frac{kp}{q}}}{\Gamma(\frac{kp}{q}+1 )} \\
	&=\sum_{k=0}^{\infty}\frac{\xi^{(\frac{k}{q}-1)p}}{\Gamma(\frac{kp}{q}-p+1 )}\\
		&=\sum_{k=0}^{q-1}\frac{\xi^{(\frac{kp}{q}-p)}}{\Gamma(\frac{kp}{q}-p+1 )}+	\sum_{k=0}^{\infty}\frac{\xi^{(\frac{kp}{q})}}{\Gamma(\frac{kp}{q}+1 )}.\\
	\end{split}
	\end{equation}
	The above equation can  further be written as
	\begin{equation}
	\frac{d^p}{d\xi^p}\mathbb{E}_{p/q}(\xi^{p/q})=\sum_{k=1}^{q-1}\frac{\xi^{-kp/q}}{\Gamma~(1-kp/q)}+\mathbb{E}_{p/q}(\xi^{p/q}).
	\end{equation}
\end{proof}

\begin{theorem}
	The function $\mathbb{E}_n(\xi^n) (n = 1, 2, 3, . . .) $ satisfies the nth order ordinary
	differential equation
	\begin{equation}
	\frac{d^n}{d\xi^n}\left( \mathbb{E}_n(\xi^n)\right) =\mathbb{E}_n(\xi^n),
	\end{equation}
	where $\xi\in \mathbb{T}.$
\end{theorem}
\begin{proof}
	For $~\upalpha>0,$
		By replacing $\xi$ by $\xi^\upalpha$ in the  equation (\ref{eq: bc m2}), we get
	\begin{equation}\label{eq:ml diff eq}
	\begin{split}
	\mathbb{E}_{\upalpha}(\xi^\upalpha)&=\sum_{k=0}^{\infty}\frac{{\xi^\upalpha}^k}{\Gamma(\upalpha k+1 )}\\ 
	&=1+\frac{\xi^\upalpha}{\Gamma (\upalpha+1)} +\frac{\xi^{2\upalpha}}{\Gamma (2\upalpha+1)}+\frac{\xi^{3\upalpha}}{\Gamma(3\upalpha+1)}+\ldots\\
	\end{split}
	\end{equation}
	By taking derivative of order $\upalpha$ on both sides of the equation (\ref{eq:ml diff eq}), 
	we get,
	\begin{equation}\label{eq:ml diff eq2}
	\begin{split}
	D^\upalpha\left(  \mathbb{E}_{\upalpha}(\xi^\upalpha)\right) &=D^\upalpha \left( 1+\frac{\xi^\upalpha}{\Gamma (\upalpha+1)} +\frac{\xi^{2\upalpha}}{\Gamma (2\upalpha+1)}+\frac{\xi^{3\upalpha}}{\Gamma(3\upalpha+1)}+\ldots\right) \\ 
	&=\frac{\Gamma (1)}{\Gamma (1-\upalpha)}\xi^{-\upalpha}+\frac{\Gamma (\upalpha+1)}{\Gamma (1)} \frac{1}{\Gamma (\upalpha+1)}+\frac{\Gamma (2\upalpha+1)}{\Gamma (\upalpha+1)} \frac{\xi^\upalpha}{\Gamma (2\upalpha+1)}\\
	&~~ + \frac{\Gamma (3\upalpha+1)}{\Gamma (2\upalpha+1)} \frac{\xi^{2\upalpha}}{\Gamma (3\upalpha+1)}+\ldots\\
	&=\frac{\Gamma (1)}{\Gamma (1-\upalpha)}\xi^{-\upalpha}+1+ \frac{\xi^\upalpha}{\Gamma (\upalpha+1)} +  \frac{\xi^{2\upalpha}}{\Gamma (2\upalpha+1)}+\ldots\\
	\end{split}
	\end{equation}
	
	Since $\frac{1}{\Gamma (1-\upalpha)}=0$ ~ for $\upalpha=n\in\mathbb{N}, $ 
	 we get from equation (\ref{eq:ml diff eq2}),
	\begin{equation}
	\begin{split}
	D^n\left(  \mathbb{E}_{n}(\xi^n)\right)&=1+ \frac{\xi^n}{\Gamma (n+1)} +  \frac{\xi^{2n}}{\Gamma (2n+1)}+\ldots\\
	&=\sum_{k=0}^{\infty}\frac{{\xi^n}^k}{\Gamma(n k+1 )}\\ 
	&= \mathbb{E}_{n}(\xi^n).
	\end{split}
	\end{equation}

\end{proof}
\section{Conclusion}
In this paper,  one parameter M-L function and its properties  in bicomplex space  has been defined from its complex counterpart. Various properties and the special cases along with recurrence relations, duplication formula, integral representation, differential relation are also derived. We intend to extend the concepts of the fractional calculus in bicomplex space using M-L function. Since bicomplex space provides a more generalized approach towards the large class of functions appearing in signal theory, electromagnetism and quantum theory.
\bibliographystyle{apalike}
	\bibliography{lref}

\end{document}